\documentclass[a4paper, 11pt]{article}
\usepackage{fancyhdr}
\usepackage{amsmath,amssymb, amsfonts}
\usepackage{enumitem}
\usepackage{mathptmx}
\usepackage[resetstyle]{phfthm}
\usepackage[urlcolor=blue, colorlinks=true, linkcolor=blue, citecolor=red]{hyperref}
\usepackage{authblk}
\theoremstyle{remark}

\DeclareMathOperator{\sech}{sech}
\newcommand{\rmd}{\mathrm}
\allowdisplaybreaks
\title{\bf An Alternate proof for a case of a Malmsten integral}
\author[ ]{Abdulhafeez Ayinde Abdulsalam}
\affil[ ]{Department of Mathematics, University of Ibadan, Ibadan, Nigeria\vspace{0.2cm}} 
\affil[ ]{\href{mailto:aabdulsalam030@stu.ui.edu.ng}{aabdulsalam030@stu.ui.edu.ng}, \href{mailto:hafeez147258369@gmail.com}{hafeez147258369@gmail.com}}

\fancypagestyle{fpg}{

\fancyhead[R]{\sf Abdulhafeez A. Abdulsalam}
\fancyhead[L]{\sf Romanian Mathematical Magazine}
\fancyhead[C]{}}
\begin{document}
\maketitle
\begin{abstract}
In this paper, a direct proof is presented for a case of a Malmsten integral. The method used in solving the integral is a direct one that the author has not come accross in any old or recent publication. Integration by parts, Laplace transform, an integral representation for the hyperbolic secant function, and the digamma representation for an alternating series are employed to derive the result.
\end{abstract}
\section{Introduction}
The integral evaluated in this paper is a case of the logarithmic integrals recently treated in Mathematical literature. Iaroslav V. Blagouchine in 2014 argued that the problem to be presented in this article was actually more older than reported (see \cite{chik}). He further clarified with sufficient and very convincing evidences that the so-called Vardi's integral (see \cite{vardi}) was actually a particular case of the considered family of integrals, first evaluated by Carl Malmsten and colleagues in 1842 (see \cite{malm}). Blagouchine in his fascinating article on the rediscovery of Malmsten's integrals presented several generalizations of the integral by contour integration method (see \cite{chik}). In this work, a direct proof is established while abstaining from methods used by Malmsten, Vardi, or Blagouchine. We hereby begin with the following proposition.
\newpage
\begin{proposition*} \normalfont\label{prop1} \leavevmode
\begin{enumerate}[label=(\alph*)]\label{prop1}
\item\label{propi} For $a \in \mathbb{R}$
$$\int_0^{\infty} \frac{\ln\left(x^2 + a^2\right)}{\cosh\left(\pi x\right)} = 2\ln\left(\frac{\sqrt{2}\Gamma\left(\frac{\left|a\right| }{2} + \frac{3}{4}\right)}{\Gamma\left(\frac{\left|a\right|}{2} + \frac{1}{4}\right)}\right).$$
\item\label{propii} $$\int_0^{\infty} \ln{x}\sech\left(x\right)\rmd{d}x = \pi\ln\left(\frac{\sqrt{2\pi}\Gamma\left(\frac{3}{4}\right)}{\Gamma\left(\frac{1}{4}\right)}\right) = \pi\ln\left(\frac{2\pi^{\frac{3}{2}}}{\Gamma\left(\frac{1}{4}\right)^2}\right).$$
\item\label{propiii} For $a \in \mathbb{R}^{+}$, $\Re(b) > 0$
$$\int_0^\infty \ln\left(ax\right) \sech\left(bx\right)\, \rmd{d}x = \frac{\pi}{b}\ln\left(\frac{2\sqrt{a}\pi^{\frac{3}{2}}}{\sqrt{b}\Gamma\left(\frac{1}{4}\right)^2}\right).$$
\end{enumerate}
\end{proposition*}

\begin{proof}\ \\
Let $\displaystyle \Delta\left(a\right) = \int_0^{\infty} \frac{\ln\left(x^2 +  a^2\right)}{\cosh\left(\pi x\right)} \, \mathrm{d}x$.\\\\
Then 
\begin{align*}
\Delta\left(a\right) &= \int_0^{\infty} \frac{\ln\left(\left|a\right| - ix\right)}{\cosh\left(\pi x\right)} \, \mathrm{d}x + \int_0^{\infty} \frac{\ln\left(\left|a\right| + ix\right)}{\cosh\left(\pi x\right)} \, \mathrm{d}x
\end{align*}
where $i =\sqrt{-1}$.
\begin{align*}
\Delta\left(a\right) &= 2\int_0^{\infty} \frac{\ln\left(\left|a\right| - ix\right)}{e^{-2 \pi x} + 1} e^{-\pi x} \, \mathrm{d}x +  2\int_0^{\infty} \frac{\ln\left(\left|a\right| + ix\right)}{e^{-2 \pi x} + 1} e^{-\pi x} \, \mathrm{d}x
\\&= -\frac{2}{\pi}\int_0^{\infty} \ln\left(\left|a\right| - ix\right) \,\, \mathrm{d}\left(\arctan\left(e^{-\pi x}\right)\right) 
\\&\qquad\qquad- \frac{2}{\pi}\int_0^{\infty} \ln\left(\left|a\right| + ix\right) \,\, \mathrm{d}\left(\arctan\left(e^{-\pi x}\right)\right)
\\&=\frac{-2i}{\pi}\int_0^{\infty} \frac{\arctan\left(e^{-\pi x}\right)}{\left|a\right| - ix}\, \mathrm{d}x + \frac{2i}{\pi}\int_0^{\infty} \frac{\arctan\left(e^{-\pi x}\right)}{\left|a\right| + ix} \, \mathrm{d}x + \ln{a}.
\end{align*}
It follows by Laplace transform that
\begin{align*}
\Delta\left(a\right) - \ln{a} &= \frac{-2i}{\pi}\int_0^{\infty} \arctan\left(e^{-\pi x}\right) \int_0^{\infty} e^{-t\left(\left|a\right| - ix\right)} \,\, \mathrm{d}t \, \mathrm{d}x 
\\&\qquad\qquad+ \frac{2i}{\pi}\int_0^{\infty} \arctan\left(e^{-\pi x}\right) \int_0^{\infty} e^{-t\left(\left|a\right| + ix\right)} \,\, \mathrm{d}t \, \mathrm{d}x
\\&= \frac{-2i}{\pi}\int_0^{\infty}e^{-\left|a\right|t}\int_0^{\infty} e^{itx} \arctan\left(e^{-\pi x}\right)\mathrm{d}x \,\, \mathrm{d}t \\&\qquad\qquad+\frac{2i}{\pi}\int_0^{\infty}e^{-\left|a\right|t}\int_0^{\infty} e^{-itx} \arctan\left(e^{-\pi x}\right)\mathrm{d}x \,\, \mathrm{d}t.
\\&\Delta\left(a\right) - \ln{a} = \frac{-2i}{\pi}\int_0^{\infty}e^{-\left|a\right|t}\int_0^{\infty} \left(e^{itx} - e^{-itx}\right) \arctan\left(e^{-\pi x}\right)\mathrm{d}x \,\, \mathrm{d}t.
\end{align*}
By Euler's formula,
$$e^{itx} - e^{-itx} = 2i\sin{\left(tx\right)}.$$
Therefore
\begin{align*}
\Delta\left(a\right) - \ln{a} &= \frac{4}{\pi}\int_0^{\infty}e^{-\left|a\right|t}\int_0^{\infty} \sin{\left(tx\right)} \arctan\left(e^{-\pi x}\right)\mathrm{d}x \,\, \mathrm{d}t 
\\&= \frac{4}{\pi}\int_0^{\infty}e^{-\left|a\right|t}\int_0^{\infty} \mathrm{d}\left(\frac{-\cos{\left(tx\right)}}{t}\right) \arctan\left(e^{-\pi x}\right) \,\, \mathrm{d}t
\\&= \frac{4}{\pi}\int_0^{\infty}e^{-\left|a\right|t}\left( \frac{\pi}{4t} - \frac{\pi}{2t}\int_0^{\infty}\frac{\cos{\left(tx\right)}}{\cosh\left(\pi x\right)} \, \mathrm{d}x\right) \,\, \mathrm{d}t 
\\&= \frac{4}{\pi}\int_0^{\infty}e^{-\left|a\right|t}\left( \frac{\pi}{4t} - \frac{1}{2t}\int_0^{\infty}\frac{\cos{\left(\frac{tx}{\pi}\right)}}{\cosh\left(x\right)}\, \mathrm{d}x\right) \,\, \mathrm{d}t
\\&= \frac{4}{\pi}\int_0^{\infty} e^{-\left|a\right|t} \left(\frac{\pi}{4t} - \frac{\pi}{4t}\sech\left(\frac{t}{2}\right)\right) \,\, \mathrm{d}t 
\\&= \int_0^{\infty} e^{-\left|a\right|t}\left(\frac{1}{t} - \frac{1}{t}\sech\left(\frac{t}{2}\right)\right) \,\, \mathrm{d}t 
\\&=  \int_0^{\infty} \left(\frac{e^{-\left|a\right|t}}{t} - \frac{2e^{-\left(\left|a\right| + \frac{1}{2}\right)t}}{t\left(1 + e^{-t}\right)} \right)\,\, \mathrm{d}t {\stackrel{\,\,t \rightarrow 2t}{=}} \int_0^{\infty} \left(\frac{e^{-2\left|a\right|t}}{t} - \frac{2e^{-\left(2\left|a\right| + 1\right)t}}{t\left(1 + e^{-2t}\right)} \right)\,\, \mathrm{d}t  
\\&{\stackrel{z \rightarrow e^{-t}}{=}} -\int_0^{1} \left(z^{2\left|a\right|} - \frac{2z^{2\left|a\right|+1}}{1 + z^2} \right)\,\, \frac{\mathrm{d}z}{z\ln{z}} = -\int_0^{1} \frac{z^{2\left|a\right|}}{\ln{z}}\left(\frac{1}{z} - \frac{2}{1 + z^2}\right)\,\, \mathrm{d}z 
\\&= -\int_0^{1} \frac{z^{2\left|a\right|}}{\ln{z}}\left(\frac{1 + z^2 - 2z}{z\left(1 + z^2\right)}\right)\,\, \mathrm{d}z = -\int_0^{1} \frac{z^{2\left|a\right|}}{\ln{z}}\left(\frac{\left(1 - z\right)^2}{z\left(1 + z^2\right)}\right)\,\, \mathrm{d}z 
\\&= \int_0^{1} \frac{z^{2\left|a\right| - 1}\left(1 - z\right)}{1 + z^2}\int_0^1 z^p \mathrm{d}p\,\, \mathrm{d}z = \int_0^{1}\int_0^1 \frac{z^{2\left|a\right| + p - 1}\left(1 - z\right)}{1 + z^2} \mathrm{d}p\,\, \mathrm{d}z  
\\&= \int_0^{1}\int_0^1 \sum_{k=0}^{\infty} \left(-1\right)^k z^{2\left|a\right| + p + 2k - 1}\left(1 - z\right) \mathrm{d}z\,\, \mathrm{d}p \\&= \int_0^{1} \sum_{k=0}^{\infty} \left(-1\right)^k \int_0^1 z^{2\left|a\right| + p + 2k - 1}\left(1 - z\right) \mathrm{d}z\,\, \mathrm{d}p 
\\&= \int_0^{1} \sum_{k=0}^{\infty} \left(-1\right)^k \left(\frac{1}{2\left|a\right| + p + 2k} - \frac{1}{2\left|a\right| + p + 2k + 1}\right)\, \mathrm{d}p 
\\&= \frac{1}{2}\int_0^{1} \sum_{k=0}^{\infty} \left(-1\right)^k \left(\frac{1}{k + \frac{2\left|a\right| + p}{2}} - \frac{1}{k + \frac{2\left|a\right| + p + 1}{2}}\right)\, \mathrm{d}p 
\\&= -\frac{1}{4}\int_0^{1} \left(\psi_0\left(\frac{2\left|a\right| + p}{4}\right) - \psi_0\left(\frac{2\left|a\right| + p}{4} + \frac{1}{2}\right) \right.
\\&\qquad\qquad\qquad\left.- \psi_0\left(\frac{2\left|a\right| + p + 1}{4}\right) + \psi_0\left(\frac{2\left|a\right| + p + 1}{4} + \frac{1}{2}\right)\right)  \mathrm{d}p
\\&= -\ln\left(\frac{\Gamma\left(\frac{2\left|a\right| + p}{4}\right)\Gamma\left(\frac{2\left|a\right| + p+ 1}{4} + \frac{1}{2}\right)}{\Gamma\left(\frac{2\left|a\right| + p}{4} + \frac{1}{2}\right)\Gamma\left(\frac{2\left|a\right| + p + 1}{4}\right)}\right)\biggr\vert_0^1 
\\&=  -\ln\left(\frac{\Gamma\left(\frac{2\left|a\right| + 1}{4}\right)\Gamma\left(\frac{2\left|a\right| + 2}{4} + \frac{1}{2}\right)}{\Gamma\left(\frac{2\left|a\right| + 1}{4} + \frac{1}{2}\right)\Gamma\left(\frac{2\left|a\right| + 2}{4}\right)}\right) + \ln\left(\frac{\Gamma\left(\frac{2\left|a\right|}{4}\right)\Gamma\left(\frac{2\left|a\right| + 1}{4} + \frac{1}{2}\right)}{\Gamma\left(\frac{2\left|a\right|}{4} + \frac{1}{2}\right)\Gamma\left(\frac{2\left|a\right| + 1}{4}\right)}\right) 
\\&= -\ln\left(\frac{\Gamma\left(\frac{2\left|a\right| + 1}{4}\right)^2\Gamma\left(\frac{2\left|a\right| + 2}{4} + \frac{1}{2}\right)}{\Gamma\left(\frac{2\left|a\right| + 1}{4} + \frac{1}{2}\right)^2 \Gamma\left(\frac{2\left|a\right| + 2}{4}\right)}\right) + \ln\left(\frac{\Gamma\left(\frac{2\left|a\right|}{4}\right)}{\Gamma\left(\frac{2\left|a\right|}{4} + \frac{1}{2}\right)}\right) 
\\&=  -\ln\left(\frac{\Gamma\left(\frac{2\left|a\right| + 1}{4}\right)^2\Gamma\left(\frac{\left|a\right|}{2} + 1 \right)}{\Gamma\left(\frac{2\left|a\right| + 1}{4} + \frac{1}{2}\right)^2 \Gamma\left(\frac{2\left|a\right| + 2}{4}\right)}\right) + \ln\left(\frac{\Gamma\left(\frac{\left|a\right|}{2}\right)}{\Gamma\left(\frac{2\left|a\right| + 2}{4}\right)}\right) 
\\&= -\ln\left(\frac{\left|a\right|\Gamma\left(\frac{2\left|a\right| + 1}{4}\right)^2}{2\Gamma\left(\frac{2\left|a\right| + 1}{4} + \frac{1}{2}\right)^2}\right) = 2\ln\left(\frac{\sqrt{2}\Gamma\left(\frac{2\left|a\right| + 1}{4} + \frac{1}{2}\right)}{\sqrt{\left|a\right|}\Gamma\left(\frac{2\left|a\right| + 1}{4}\right)}\right).
\end{align*}
Hence
\begin{align}
\Delta\left(a\right) &= 2\ln\left(\frac{\sqrt{2}\Gamma\left(\frac{2\left|a\right| + 1}{4} + \frac{1}{2}\right)}{\sqrt{\left|a\right|}\Gamma\left(\frac{2\left|a\right| + 1}{4}\right)}\right) + \ln{\left|a\right|} = 2\ln\left(\frac{\sqrt{2}\Gamma\left(\frac{\left|a\right|}{2} + \frac{3}{4}\right)}{\Gamma\left(\frac{\left|a\right|}{2} + \frac{1}{4}\right)}\right),\label{rzq}
\end{align}
thereby proving \ref{propi} in \autoref{prop1}.
To prove \ref{propii} in \autoref{prop1}, take the limit of \eqref{rzq} as $a \to 0$.\\\\
Thus
$$\int_0^{\infty} \frac{\ln\left(x\right)}{\cosh\left(\pi x\right)} \, \mathrm{d}x = \ln\left(\frac{\sqrt{2}\Gamma\left(\frac{3}{4}\right)}{\Gamma\left(\frac{1}{4}\right)}\right),$$
which implies
$$ \int_0^{\infty} \ln{x}\sech\left(\pi x\right)\mathrm{d}x= \ln\left(\frac{\sqrt{2}\Gamma\left(\frac{3}{4}\right)}{\Gamma\left(\frac{1}{4}\right)}\right)$$
$$ \frac{1}{\pi}\int_0^{\infty} \ln\left(\frac{x}{\pi}\right)\sech\left(x\right)\mathrm{d}x = \ln\left(\frac{\sqrt{2}\Gamma\left(\frac{3}{4}\right)}{\Gamma\left(\frac{1}{4}\right)}\right)$$
$$\int_0^{\infty} \ln\left(\frac{x}{\pi}\right)\sech\left(x\right)\mathrm{d}x = \pi\ln\left(\frac{\sqrt{2}\Gamma\left(\frac{3}{4}\right)}{\Gamma\left(\frac{1}{4}\right)}\right)$$
$$\int_0^{\infty} \ln{x}\sech\left(x\right)\mathrm{d}x - \ln{\pi}\int_0^{\infty} \sech\left(x\right)\mathrm{d}x = \pi\ln\left(\frac{\sqrt{2}\Gamma\left(\frac{3}{4}\right)}{\Gamma\left(\frac{1}{4}\right)}\right)$$
$$\int_0^{\infty} \ln{x}\sech\left(x\right)\mathrm{d}x - \frac{\pi}{2}\ln{\pi} = \pi\ln\left(\frac{\sqrt{2}\Gamma\left(\frac{3}{4}\right)}{\Gamma\left(\frac{1}{4}\right)}\right)$$
$$\int_0^{\infty} \ln{x}\sech\left(x\right)\mathrm{d}x= \pi\ln\left(\frac{\sqrt{2\pi}\Gamma\left(\frac{3}{4}\right)}{\Gamma\left(\frac{1}{4}\right)}\right)=  \pi\ln\left(\frac{2\pi^{\frac{3}{2}}}{\Gamma\left(\frac{1}{4}\right)^2}\right).$$
To prove \ref{propiii} in \autoref{prop1},
\begin{align*}
\int_0^\infty \ln\left(ax\right) \sech\left(bx\right)\, \rmd{d}x &= \int_0^\infty \ln\left(ax\right) \sech\left(bx\right)\, \rmd{d}x
\\&= \frac{1}{b}\int_0^\infty \ln\left(\frac{au}{b}\right)\sech\left(u\right) \, \rmd{d}x
\\&=  \frac{1}{b}\int_0^\infty \ln\left(u\right)\sech\left(u\right) \, \rmd{d}x - \frac{\ln{\left(\frac{b}{a}\right)}}{b}\int_0^\infty \sech\left(u\right) \, \rmd{d}x
\\&= \frac{\pi}{b}\ln\left(\frac{2\pi^{\frac{3}{2}}}{\Gamma\left(\frac{1}{4}\right)^2}\right) - \frac{\ln{\left(\frac{b}{a}\right)}}{b}\cdot\frac{\pi}{2}
\\&= \frac{\pi}{b}\ln\left(\frac{2\sqrt{a}\pi^{\frac{3}{2}}}{\sqrt{b}\Gamma\left(\frac{1}{4}\right)^2}\right).
\end{align*}
\end{proof}

\bibliographystyle{unsrt}
\bibliography{RMM}

\end{document}